%
%
\input amstex
\documentstyle{jams}
\NoBlackBoxes

\topmatter
\title Oort's conjecture for  $A_g \otimes {\Bbb C}$ \endtitle
\author Sean Keel and Lorenzo Sadun \endauthor
\leftheadtext{Sean Keel and Lorenzo Sadun}%
\rightheadtext{Oort's conjecture for $A_g \otimes {\Bbb C}$}%

\address Department of Mathematics, University of Texas at Austin, 
Austin, Texas, 78712 \endaddress


\email keel\@math.utexas.edu\endemail

\thanks The first author was partially supported by
NSF grant DMS-9988874\endthanks 

\address Department of Mathematics, University of Texas at Austin, 
Austin, Texas, 78712 \endaddress

\email sadun\@math.utexas.edu\endemail

\thanks The second author was partially supported by
Texas ARP grant 003658-152. \endthanks

\issueinfo{00}
  {0}
  {Xxxx}
  {0000}

\subjclass Primary 14K10 \endsubjclass



\abstract
We prove the conjecture of Oort that a compact subvariety
of the moduli space of principally polarized Abelian
varieties of genus $g$ has codimension strictly greater than 
$g$, in characteristic zero, for $g \geq 3$.
\endabstract
\date January 17, 2003 \enddate
\endtopmatter

\document

\magnification\magstep1
\NoRunningHeads
\pageheight{9 truein}
\pagewidth{6.3 truein}
\baselineskip=15pt
\catcode`@=11
\def\logo@{\relax}
\catcode`@=\active
\def\bP{\Bbb P}

\def\cc{\tilde{c}}

\def\bH{\Bbb H}
\def\bC{\Bbb C}
\def\bQ{\Bbb Q}

\def\bE{\Bbb E}

\def\cR{\Cal R}
\def\ring#1.{\Cal O_{#1}}
\def\Sp{\operatorname{Sp}}
\def\U{\operatorname{U}}

\def\ag{A_g}

\def\dr{\operatorname{DR}}
\def\mg{\overline{M}_g}
\def\mgc{M_g^c}
\def\omg{M_g}
\def\mgco{M_{g,1}^c}
\def\omgo{M_{g,1}}
\def\ch{\operatorname{CH}}

\def\hg{{\Cal H}_g}
\def\vhg#1.{{\Cal H}_{#1}}

\def\im{\operatorname{Im}}
\def\re{\operatorname{Re}}

\def\bar#1.{\overline{#1}}

\def\Hom{\operatorname{Hom}}
\def\cX{\Cal X}

\def\barW{\bar {W}.}

\def\bF{\Bbb F}
\def\bR{\Bbb R}
\def\dim{\operatorname{dim}}
\def\codim{\operatorname{codim}}
\def\det{\operatorname{det}}
\def\bP{\Bbb P}

\def\ker{\operatorname{ker}}
\def\im{\operatorname{im}}
\def\Hom{\operatorname{Hom}}

\def\bZ{\Bbb Z}
\def\Gal{\operatorname{Gal}}
\def\lag{\goth g}
\def\Sp{\operatorname{Sp}}
\def\sp{\operatorname{sp}}

\def\U{\operatorname{U}}
\def\u{\operatorname{u}}

\def\ag{A_g}

\subhead \S 1. Introduction and Results \endsubhead
Let $\ag$ be the moduli space of principally polarized
Abelian varieties and $\bE$ be the Hodge bundle. 

\proclaim{1.1 Main Theorem} Let $X$ be a compact
subvariety of $\ag \otimes \bC$ such
that $c_i(\bE)|_X =0$ in $H^{2i}(X,\bR)$
for some $g \geq i \geq 0$. 
Then 
$$
\dim(X) \leq \frac{i(i-1)}{2}, \leqno(1)
$$
with strict inequality if $i \geq 3$.
\endproclaim

The
Grothendieck Riemann Roch theorem
implies that $c_g(\bE)$ is trivial
on $\ag$, see \cite{G99,2.2}, thus
we have the following corollary, conjectured
by Oort, \cite{OG,3.5}:

\proclaim{1.2 Corollary} There is no compact
codimension $g$ subvariety of $\ag \otimes \bC$,
for $g \geq 3$.
\endproclaim

Corollary 1.2 is striking in that it fails in positive
characteristic: Oort has shown that the
locus $Z \subset \ag \otimes \bF_p$ of Abelian
varieties of $p$-rank zero is complete and pure
codimension $g$. 

We also obtain:
\proclaim{1.2.1 Corollary} There is no compact codimension
$g$ subvariety of $\mgc \otimes \bC$ for $g \geq 3$.
\endproclaim
This  is a formal consequence of (the $g=3$ case
of) (1.2), and the bounds of \cite{Diaz87} on
dimensions of compact subvarieties of $\mg$, see
\S 8.

Inequality (1) has been obtained previously,
in all characteristics, by van der Geer, \cite{G99}. The
improvement to strict inequality may seem
rather humble, but in relation to Faber's conjectures, 
it is quite significant:

Faber has made the surprising conjecture that the subring
$\cR^*(\omg)\subset \ch^*(\omg)_{\bQ} \otimes \bQ$ generated by
tautological classes {\it looks like} the rational cohomology of a
smooth projective variety of dimension $g-2$, e.g., it satisfies
Poincare duality and Grothendieck's standard conjectures. See
\cite{Faber99} for the precise statement, and very compelling
evidence. One naturally wonders if there is a projective $g-2$
dimensional variety with the right cohomology.  The form of the
conjectures, specifically his evaluation maps, \cite{FP00,0.4},
strongly suggest looking for a compact codimension $2g-1$ subvariety
of $\omg$, whose cohomology class in $\mg$ is a multiple of
$c_g(\bE)c_{g-1}(\bE)$. There are analogous speculations for $\omg^c$
(stable curves of compact type), \cite{FP00,0.5}, that lead one to
hope for a compact subvariety of codim $g$, with cohomology class a
multiple of $c_g(\bE)$. Since $\bE$ comes from $\ag$, it is natural to
hope that the subvarieties do as well. More compellingly, there are
natural analogs of Faber's conjectures for $\ag$, and $\ag^{\circ}$
(the open subset of $\ag$ whose complement correspond to products of
lower dimensional polarized varieties -- note the image of the Torelli
map is is a closed subset of $\ag^{\circ}$), and corresponding
evaluation maps which themselves raise the hope of compact
subvarieties, of codimension $g$ in $\ag$, and of codimension $2g-1$ in
$\ag^{\circ}$. Indeed the tautological ring $\cR^*(\ag)$ has been
computed, \cite{G99} and \cite{EV02}, and the analog of Faber's
conjecture holds, so if in particular one has a compact subvariety $Z
\subset \ag$ with the desired cohomology class, then Poincare duality
already gives that $\cR^*(\ag) \rightarrow H^*(Z,\bQ)$ is injective.

Unfortunately, by Corollary 1.2 the hoped for subvariety of $\ag$ does
not exist in characteristic zero. Theorem 1.1 suggests a similar result may
hold for $\ag^{\circ}$: It is natural to wonder if $c_{g-1}(\bE)$
vanishes on $\ag^{\circ}$.  If it does, then by Theorem 1.1, with $i = g-1$,
there is no compact subvariety of $\ag^{\circ} \otimes \bC$ of
codimension $2g-1$, for $g \geq 4$.

As noted, Oort's example $Z \subset \ag \otimes \bF_p$ 
violates the analog of 
(1.2) for characteristic $p$. As Oort pointed out
to us, since it is known that $Z$ does not lift
to char $0$, uniqueness of $Z$ would give a different,
purely char $p$ proof of (1.2). In this sense (1.2) 
is evidence for this uniqueness, which would obviously
be attractive from the point of view of Faber's 
conjecture. In \S 7 we point out an intriguing parallel
between our proof in char $0$ and Oort's example; perhaps
it can be exploited. 

Finally, Theorem 1.1 is sharp for $i=2$ and $i=3$.
\def\vag#1.{A_{#1}}
There exists a compact
surface $Z \subset \vag 3.$ and a compact curve 
$C \subset \vag 2.$. Taking products with
fixed elliptic curves one then obtains
a compact surface $Z \subset \ag$, $g \geq 3$ and a 
compact curve $C \subset \ag$, $g \geq 2$ with 
$c_3(\bE)|_Z =0$, $c_2(\bE)|_C =0$. 

After submitting this paper for publication 
we learned from
E. Colombo that the $i=g$ case of Corollary (2.6),
the crucial technical step in the proof of Corollary (1.2), 
has been previously obtained by E. Izadi, \cite{Izadi98},
using methods of Colombo and Pirola, \cite{CP90}. 
The main argument of Colombo-Pirola and Izadi 
is entirely
different from ours (and quite ingenious). It is more
elementary (no use is made of the vanishing of
$c_g(\bE)$) and also more algebaic, so might be
useful for showing uniqueness of 
Oort's $Z \subset \ag \otimes \bF_p$.

Thanks: We would like to thank 
D.~Freed, 
J.~Tate, F.~Villegas, G.~van der Geer,
R.~Hain, C.~Faber, and
C. Teleman.
F.~Voloch helped us with the
proof of (7.1), and M.~Stern sent us copious
explanations of the Satake compactification. 
A.~Reid, and D. ~Witte gave us lots of help
with discrete groups. Finally the main ideas
in the proof of (6.3) are due to Scot Adams. 

\subhead \S 2. Proof of Main Theorem \endsubhead

In this section we give the main line of our argument for
Theorem 1.1. The various constituent results will be
proved in subsequent sections. 

We use throughout the (orbifold) cover
$q: \hg \rightarrow \ag$ by Siegel's generalized
upper half space. (We could as well replace $\ag$ by
a finite branched cover, and $q$ by an honest 
(unbranched) cover. Above and throughout
the paper any statement
about $\ag$ should be interpreted in the
orbifold sense.)
$\hg$ is by definition the
space of symmetric 
complex $g \times g$ 
matrices with positive
definite imaginary part, while $S_g$ denotes the space of all 
symmetric $g \times g$
complex matrices. 

\definition{2.1. Chern Forms} $\bE$ has a tautological Hermitian
metric: A point of $\ag$ is given by a pair $(V/L,<,>)$, where $L
\subset V$ is a lattice in a $g$-dim complex vector space, and $<,>$
is a positive definite Hermitian form (such that $\im <,>|_L$ is
integral and unimodular). $\bE$ is the Hermitian bundle with fibre
$V^*$ and metric dual to $<,>$. (We abuse notation and denote this
dual metric, and any other metric induced from $<,>$, by
$<,>$). Griffiths showed that the curvature of $(\bE,<,>)$ is non-negative,
\cite{Griffiths84}.  Our proof is based on the tension between this
non-negativity and the vanishing of $c_g(\bE)$ noted above: As $\bE
\geq 0$ it's natural to hope that $c_k(\bE)$ is represented by a
non-negative form $\cc_k$.  (Recall a real $(k,k)$ form is called non-negative
if its restriction to every complex $k$-plane in the tangent bundle is
a non-negative multiple of the volume form.)  As we indicate in a
moment, this is indeed the case. Thus if $Y \subset \ag$ is a compact
$k$-fold and
$$
0 = c_k(\bE) \cap [Y] = \int_Y \cc_k, \leqno(2)
$$
then $\cc_k|_Y$ vanishes identically
(on the smooth locus), and we get pointwise
information on the tangent space of $Y$.
\enddefinition

Here is how we obtain $\cc_k$.
By Hodge theory,
$\bE$ sits in an exact sequence
$$
0 \rightarrow \bE \rightarrow \bH \rightarrow \bE^* 
\rightarrow 0, \leqno(3)
$$
where $\bH$ is a flat bundle (with fibre $H^1(A,\bC)$ at
the point $[A] \in \ag$ given by the Abelian variety
$A$), whence the equality of cohomology classes
$$
c(\bE) \cdot c(\bE^*) =1. \leqno(4) 
$$
Thus $c(\bE)$ is represented by the inverse 
Chern form of $(\bE^*,<,>)$. This is computed as follows:

\definition{2.2 Segre forms} Let $(F,<,>)$ be a rank $g$
Hermitian bundle, on a complex manifold $B$, with
curvature $\Omega$, let $G = \frac{1}{2\pi i} \Omega$, and let 
$p: \bP(F) \rightarrow B$ be the projective
bundle of lines in $F$. Let
$$
A: T_{\bP(F)} \rightarrow T_{\bP(F)/B} \leqno(5)
$$
(where $T_{\bP(F)/B}$ is the vertical tangent
space, i.e., the kernel of $dp$) be the
second fundamental form for the tautological
sublinebundle
$\ring.(-1) \subset p^*(F)$. Let
$s(F,<,>)= \sum s_k(F,<,>)$ denote the Segre from, i.e.,
$$
s_k(F) = p_*(c_1(\ring.(1),<,>)^{g-k+1}). \leqno(6)
$$

\enddefinition

\proclaim{2.3 Theorem} 
The kernel of $A$ 
defines a horizontal subspace, complementary
to the vertical tangent space, which is
identified with $p^*(T_{B})$ by $dp$. Under
this splitting of $T_{\bP(F)}$ the first Chern form of
$\ring.(1)$ at a line $L \subset F_b$, 
spanned by a unit vector $v$, is 
$$
FS + <Gv,v>, \leqno(7) 
$$
where $FS$ is the Fubini-Study form
(on the vertical tangent space, normalized to represent the hyperplace class).
The Segre form $s_k(F,<,>)$ equals
$$
{g+k-1 \choose k} 
\int_{\bP(F_b)} \left ({<G v,v> \over <v,v>}\right )^k dVol_{\bP(F_b)}, 
\leqno(8)
$$
and the total Segre form $s(F,<,>)$ equals  $c(F,<,>)^{-1}$ (pointwise!) 
in the algebra of even degree forms.
\endproclaim

It is well known that the Segre class represents the
inverse Chern class in cohomology, see e.g. \cite{BT82,pg 269},
\cite{Fulton80,3.1}.
We do not know if this equality at
the level of forms, or 
the rather canonical expressions above, have been previously
observed. 

The tangent space $T_{p} \ag$ is canonically
identified with $S(\bE_p,\bE_p^*)$, the
space of symmetric linear maps from 
$\bE_p$ to its dual. 
Theorem 2.3 applied to $(\bE^*,<,>)$, together
with Griffiths' formula for the curvature, 
imply:
\proclaim{2.4 Corollary} 
The $k^{th}$ Chern class of $\bE$ is represented
by $\tilde c_k := s_k(\bE^*,<,>)$. 
The form $\tilde c_k$ is non-negative and 
vanishes on a (complex) $i$-plane $Y$ in
$T_p \ag = S(\bE_p,\bE_p^*)$ if and only if, for every $v \in \bE_p$, 
the evaluation map $e_v: Y \to \bE^*$ fails to be injective. 

For each line $L \subset \bE_p$, choose 
$0 \neq w \in L$. 
Let $<,>_L$ on $S(\bE_p,\bE_p^*)$ be 
the Hermitian form $e_{w}^*(< , >_{\bE_p^*})/<w,w>$, i.e.,
$$
<a,b>_L = { <e_{w}(a),e_{w}(b)> \over <w,w>}. \leqno(9)
$$
Up to a positive constant, 
the restriction of $\tilde c_k$
to $p$ is the average of the $k^{th}$
wedge powers of $-Im(< ,>_L)$ over all lines
$L \subset \bE_p$. \endproclaim

We prove (2.3-4) in \S 3. 

Given compact $X \subset \ag$ with $c_i(\bE)|_X =0$,
we may apply Corollary 2.4 to every $i$-plane in the tangent
space at smooth points of $X$. (Since $\ag$ is quasi-projective, each
$i$-plane is realized as a tangent plane to
a compact $i$-dimensional subvariety.)

We reformulate this using the following linear
algebra result, proved in \S 4.
\proclaim{2.5 Theorem} 
Let $V$ be 
a complex vector space of dimension $g \geq 1.$ Let
$i$ be an integer $g \geq i \geq 0$. 
Let $X \subset S(V,V^*)$ be a linear subspace
of dimension at least 
$\max( i, \frac{i(i-1)}{2})$ 
such that for each 
$v \in V$ and every $i$-dimensional subspace
$Y \subset X$ the evaluation map
$e_v:Y \rightarrow V^*$ fails to be injective. 

Then $i \geq 3$, $X = W^{\perp}$ for some $g-i+1$-dimensional
linear subspace $W \subset V$,
and $X$ has dimension exactly
$\frac{i(i-1)}{2}$. \endproclaim

Here $W^{\perp}\subset S(V,V^*)$ is the set of maps whose 
kernels contain $W$.

Combining this with Corollary 2.4 we find:

\proclaim{2.6 Corollary} Let $g \geq i \geq 3$. 
Let $X \subset \ag$
be a compact subvariety of dimension
at least $\frac{i(i-1)}{2}$ such that $c_i(\bE)|_X =0$.
Then $X$ has dimension exactly $\frac{i(i-1)}{2}$. 
Furthermore, 
at each smooth point $p \in X$,
$T_p X = W^{\perp} \subset S(\bE_p,\bE_p^*)$
for some $g - i + 1$ dimensional subspace
$W(p) \subset \bE_p$. 
\endproclaim

Let $X \subset \ag$ be as in Corollary 2.6, and let 
$\cX$ be an irreducible component of
$q^{-1}(X)$. By Corollary 2.6,  
$$
T_A X = W(A)^{\perp} \subset S_g \leqno(10)
$$
for each smooth $A \in X$. 

In \S 5 we prove:
\proclaim{2.7 Theorem}
The subspace $W(A)$ is constant (and henceforth denoted $W$).  
If $\cX$ is an irreducible component of $q^{-1}(X)$,
then for some fixed $\tau \in \hg$,
$\cX \subset \hg$ is the affine space
$$
X(W,\tau) :=\{ M \in \hg| M |_W = \tau|_W \}. \leqno(11)
$$ 
\endproclaim

Finally, in \S 6 we show:

\proclaim{2.8 Theorem}
The image of
$X(W,\tau)$ in $\ag$ cannot be compact.
\endproclaim

\subhead \S 3 Segre Forms \endsubhead
\demo{Proof of Theorem 2.3}
The restriction of $A$ to the vertical tangent
space is the identity, so $A$ is indeed surjective.
Now equation (7) follows by the definition of $A$, see
\cite{GH78,pg 78}.
What remains is to show that the pointwise inverse of the total Chern form 
$c(F,<,>)$ is given by the formula (8), and equals the
Segre form.  

If $M$ is a positive 
Hermitian form on a $g$-dimensional complex vector space $V$ with inner
product $<,>$, then
$$ \pi^{-g} \det(M) \int_{V} e^{-<Mv,v>} dVol_V = 1. \leqno(12) $$
This identity remains true if $M$ takes values in the commutative 
ring of even degree forms at a point, as long as the scalar part of 
$M$ remains positive-definite.

Now let $V = F_p$, and let $M = I -G$. Then
$$ \eqalign{
 [c(F,<,>)]^{-1} = & \pi^{-g} \int_{F_p} e^{-<(I-G)v,v>}  dVol_{F_p}\cr
= & \pi^{-g}  \int_{F_p} \sum_k {<Gv,v>^k \over k!}  
e^{-<v,v>} dVol_{F_p}\cr
(13)\phantom{MMMMMMM} = & \sum_k {g+k-1 \choose k} {1 \over Vol(S^{2g-1})} 
\int_{S^{2g-1}} <Gv,v>^k d Vol_{S^{2g-1}} 
\cr
= & \sum_k  {g+k-1 \choose k} {1 \over Vol(\bP(F_p))} 
\int_{\bP(F_p)} \left ({<Gv,v> \over <v,v>}\right )^k dVol_{\bP(F_p)}
\cr
= &
\sum_k \int_{\bP(F_p)} \left (FS + {<Gv,v>\over <v,v>}\right )^{g+k-1}.} 
$$
In going from the second line to the third, we integrate over radius, using 
the facts that $\int_0^\infty r^{2k} e^{-r^2} r^{2g-1}dr = (k+g-1)!/2$, 
and that the volume of $S^{2g-1}$ is $2\pi^g/(g-1)!$.

\qed \enddemo

\demo{Proof of Corollary 2.4}
The second paragraph implies the first, since
$< ,>_L$ is semi-positive, and its restriction
to an $i$-plane $Y$  fails to be positive
iff $e_{v^*}|_Y$ fails to be injective. Thus
${\tilde c_k}|_Y$ is non-negative
and vanishes iff $e_{v^*}$ fails to be injective for all $v$.

What remains, then, is to compute the curvature of $\bE^*$, from which 
we obtain the Segre form $s_k$ by equation (8). This is given by
Griffiths' formula \cite{Griffiths84}, but can also be computed from a
direct calculation, which we include for the reader's convenience:
The metric on $\bE^*$ is given by the inner product 
$h=(\tau_I)^{-1}$, 
where $\tau_I$ denotes the imaginary part of $\tau$,
see \cite{Kempf91,pg 59}. Thus the curvature is given by:
$$ 
\eqalign{ \Omega = &
\bar {\partial}. (h^{-1} \partial h )\cr
= & -\bar \partial. [(\partial \tau_I) \tau_I^{-1}] \cr
= & - (\partial \tau_I) \tau_I^{-1} 
(\overline{\partial \tau_I})\tau_I^{-1}
\cr
= & -\frac{1}{4} (\partial \tau) \tau_I^{-1} 
(\overline{\partial \tau})
\tau_I^{-1}.}\tag{14}
$$
If $v \in \bE^*$, then $w:= \tau_I^{-1} \bar v. = <\cdot,v> \in \bE$, and
$$\eqalign{<Gv,v> = & {i \over 8\pi} \bar v^t. \tau_I^{-1} (\partial \tau)
\tau_I^{-1} (\overline{\partial \tau})
\tau_I^{-1} v \cr
= & {i \over 8\pi} < e_{w}(\partial \tau), e_{w}(\partial \tau)>   \cr
= & {-1 \over 4 \pi} Im \left ( e_w^* < , >_{\bE_p^*}\right ).} \tag{15}
$$
Note that $<v,v>=<w,w>$, and the unit sphere in $\bE^*_p$ is
naturally identified with the unit sphere in $\bE_p$. As a result, averaging
$<Gv,v>^k/<v,v>^k$ over all values of $v$
is equivalent to averaging $(-Im(<,>_L)/4\pi)^k$ over all lines 
in $\bE_p$. \qed
\enddemo

Remark: It is natural to wonder whether or not 
the Chern form
$c_k(\bE,<,>)$ and the Segre form
$s_k(\bE^*,<,>)$ are actually equal. We checked
that this holds for $k =1,2$, and suspect that it holds in
general. 
\subhead \S 4 Linear algebra \endsubhead

As before, let $S(V,V^*)$ be the space of 
symmetric maps from the vector space $V$
to its dual. For each $v \in V$, let $v^\perp \subset S(W,W^*)$ 
be the set of maps that vanish on $v$. Let 
$R_k \subset S(W,W^*)$ be the
locus of maps whose image has dimension
at most $k$. (This is a closed subset,
not a vector subspace.)

We prove Theorem 2.5 in stages, beginning with the following standard fact:

\proclaim{4.1 Lemma} Let $M\in R_k \setminus R_{k-1}$.
$R_k$ is smooth at $M$ with tangent space
$$
T_M R_k = \{N \in S(V,V^*)| N(\ker(M)) \subset \im(M) \}. \leqno(16)
$$
\endproclaim

\proclaim{4.2 Corollary} 
Let $Y \subset S(V,V^*)$ be a $2$-plane spanned
by rank $1$-maps. Then $Y \subset R_2$, but $Y \not \subset R_1$.
\endproclaim
\demo{Proof} Let $M,N \in R_1$ be a basis of $Y$, and let $T \in Y$.
Then $\ker(M) \cap \ker(N) \subset \ker(T)$, so 
$T \in R_2$. If $Y \subset R_1$, however, then 
$M \in T_N R_k$, so $M(\ker N) \subset \im(N)$. Using
that they are each rank one, this implies that $M$ and $N$
are dependent, a contradiction. \qed \enddemo

\proclaim{4.3 Lemma} Let $Y \subset S(V,V^*)$ 
be a linear subspace of $R_2$. 
Let $M \in Y$, $v \in V$ be such that $M(v)(v) \neq 0$.
Then $v^{\perp} \cap Y \subset R_1$ and has dimension
at most one.
\endproclaim
\demo{Proof}
Suppose $N \in v^{\perp}\cap Y$ has rank $2$. 
Since $M \in T_N R_2$, Lemma 4.1 implies that 
$M(v) = N(w)$ for some $w \in V$.
But then $M(v)(v) = N(v)(w) =0$, a contradiction.
Thus $v^\perp \cap Y \subset R_1$. Now apply Corollary 4.2 to $v^\perp
\cap Y$. 
\qed \enddemo

It will be convenient to prove the cases $i=1,2$ of Theorem 2.5
separately:
\proclaim{4.4 Proposition} If $i<3$, then the
hypotheses of Theorem 2.5 are never met.
\endproclaim
\demo{Proof} We consider an $i$-dimensional subspace $Y$ on which, 
for all nonzero $v$, 
$e_v$ fails to be injective. That is, every vector in $V$ is in the kernel of 
some nonzero element of $Y$, and we have
$$
V = \bigcup_{Q \in \bP(Y)} \ker(Q). \leqno(17)
$$
If $i=1$, this means that all of $V$ is in the kernel of a nonzero
matrix, clearly impossible.  If $i=2$, this means that a
$g$-dimensional vector space is the union of a 1-dimensional family of
proper subspaces.  Thus almost every element of $Y$ has a
$(g-1)$-dimensional kernel, i.e., is of rank 1.  But this contradicts
Corollary 4.2. \qed\enddemo

\proclaim{4.5 Proposition} Let $i\ge 3$. If Theorem 2.5 holds for all 
spaces $X$ of dimension equal to $i(i-1)/2$, then the hypotheses of
the theorem cannot be satisfied by any $X$ of
dimension greater than $i(i-1)/2$. 
\endproclaim

\demo{Proof} Let $X$ have dimension greater than $i(i-1)/2$, satisfying
the hypotheses of Theorem 2.5, and let
$X'$ be a subspace of dimension $i(i-1)/2$.  The hypotheses of the theorem
then apply to $X'$, so $X' \subset W^\perp$ for some $g-i+1$ dimensional 
subspace $W \subset V$. Since their dimensions are equal, $X' = W^\perp$,
so we can find $i-1$ rank-1 elements $x_1,\ldots x_{i-1}$ of $X'$ whose 
ranges are linearly independent. Let $x_i \in X - X'$, and let $Y$ be
the span of $x_1,\ldots, x_n$.  It is then straightforward to find a vector
$v \in V$ for which $e_v$ is injective on $Y$, which is a contradiction.  
\enddemo

In light of Proposition 4.5, we will henceforth only consider spaces $X$
of dimension equal to $i(i-1)/2$. 

\proclaim{4.6 Proposition} Theorem 2.5 holds for $i=3$.
\endproclaim

\demo{Proof} Since $i=i(i-1)/2=3$, the only $i$-dimensional subspace
$Y$ of $X$ is $X$ itself.  We first show that $Y \subset R_2$. This
follows from equation (17) by counting dimensions. If $Y \not \subset
R_2$, then there must exist a codimension-one locus $Z' \subset
\bP(Y)$ of rank $1$ quadrics, and correspondingly a codim one closed
cone $Z \subset Y \cap R_1$.  But then, by Corollary 4.2, 
$Y \subset R_2$, which is a contradiction.

Having shown that $Y \subset R_2$, we use Lemma 4.3 to construct
rank-1 elements of $Y$.  Since $Y$ is nonzero, we can find $M_0 \in Y$
and $v_0 \in V$ such that $M_0(v_0)(v_0) \ne 0$.  Since $e_{v_0}$ in 
not injective, $v_0^\perp \cap Y$ is nonempty, and contains a rank-1
element $M_1$.  Let $v_1 \in V$ be such that $M_1(v_1)(v_1) \ne 0$. 
By Lemma 4.3, we can then find a rank-1 element element $M_2 \in v_1^\perp$.

Clearly $M_1$ and $M_2$ are independent.  Let $Z$ be the $2$-plane
they span, and let $W = \ker(M_1) \cap \ker(M_2)$.  $W \subset V$ is
codimension two and $Z \subset W^{\perp}$. Let $H$ be another element
of $Y$, independent of $M_1$ and $M_2$, so that $Y$ is the span of
$M_1$, $M_2$, and $H$.

By Proposition 4.4, there
exists $t \in V$ such that
$$
e_t: Z \rightarrow V^*  \leqno(18)
$$
is injective.  Injectivity is an open condition, so $e_t$ is injective
on $Z$ for general $t$.  However, by assumption $e_t$ is {\it never} 
injective on $Y$. 
Thus for each $t$ there exists $M_3 \in Z$ such that 
$M_3(t)=H(t)$. However, that implies that, for every $w \in W$, 
$$ H(w)(t) = H(t)(w) = M_3(t)(w)=M_3(w)(t) = 0. \leqno(19) $$
Since this is true for general $t$, $H(w)$ must be zero, so $H \in W^\perp$.
But then $Y \subset W^\perp$.

Finally, the dimension of $W^\perp$ is exactly 3, so $X= Y = W^\perp$. 
\qed \enddemo




For $i>3$, we will prove Theorem 2.5 by induction on $i$.  
The key to the inductive step is the following lemma:
 
\proclaim{4.7 Lemma} Suppose that Theorem 2.5 applies for all values of
$i \le k$, and that $X$ is a subspace of $Sym(V,V^*)$ satisfying the
hypotheses of theorem (1.4) with $i=k+1$. Let $M \in X, v \in V$ be
such that 
$M(v)(v) \neq 0$. There is a $g-i+2$-dimensional
subspace $W \subset V$, containing $v$, such that
$$
W^{\perp}  
= v^{\perp} \cap X \subset S(V,V^*) \leqno(20)
$$
In particular $v^{\perp} \cap X$ contains non-zero elements
of rank one. \endproclaim
\demo{Proof} 
If $W \subset V$ is a linear subspace,
and $W^{\perp} \subset S(V,V^*)$ is
the subspace of maps $f$ with $f|_W =0$, then
the natural map $W^\perp \rightarrow \Hom(V/W,(V/W)^*)$
identifies $W^\perp$ with $S(V/W,(V/W)^*)$.

Let $X_v:= v^{\perp} \cap X$. 
By assumption, $e_v$ is not injective on every $i$-dimensional subspace
$Y \subset X$.  This implies that the rank of $e_v : X \to V^*$ is at most
$i-1$, so the kernel of $e_v$ has codimension at most $i-1$. Thus 
$$
\align
\dim(X_v) = \dim(X) - \codim(X_v,X) 
&\geq \frac{i(i-1)}{2} - (i-1) \\
&=\frac{(i-1)(i-2)}{2}. \tag{21}
\endalign 
$$

We claim that 
$X_v \subset S(V/v,(V/v)^*)$ 
satisfies the evaluation condition of Theorem 2.5. 
If not, then for some $w \neq 0$, the range of 
$e_w$ on $X_v$ has dimension $i-1$, and so equals the
range of $e_w$ on all of $Y$. 
Note that $e_w = e_{w + v}$ on $X_v \subset v^\perp$,
so (after possibly replacing $w$ by $w + v$) we
can assume $M(w)(v) \neq 0$. But now
$M(w) = N(w)$ for some $N \in v^{\perp}$, thus
$M(w)(v) = N(w)(v)=N(v)(w)= 0$, a contradiction. 

We now apply Theorem 2.5 to $X_v$. Viewed as an element of $S(V/v, (V/v)^*)$,
$X_v$ is the orthogonal complement of a $(g-i+1)$-dimensional subspace
$W_0 \subset V/v$.  This is tantamount to equation (20), where
$W$ is the preimage of $W_0$.
\qed \enddemo

\demo{Proof of Theorem 2.5} In light of Propositions 4.4, 4.5 and 4.6,
we need only prove Theorem 2.5 for $i>3$ and for $X$ of dimension
exactly $i(i-1)/2$.  We do this by induction on $i$, relying on the previously
proven base case $i=3$.

To establish the induction, we apply Lemma 4.7 to get a rank one element
$M \in X$. 
Choose $v \in V$ with $M(v)(v) \neq 0$, and apply Lemma 4.7 again, so
that
$$                
\barW^{\perp} = X_v \subset X \leqno(22)
$$
for some $g-i+2$-plane $\barW \subset V$ that contains $v$.
Now let $W = \ker(M) \cap \barW$.  We claim that $X = W^\perp$.

To see this, let $Z \subset X$ be the span of $\barW^\perp$ and $M$,
and let $t$ be an element of $V$ that is neither in $\barW$ nor in
the kernel of $M$. (This is an open condition).  Note that $M(t)$ is a
nonzero multiple of $M(v)$, so $M(t)(v) \ne 0$.

Since $t \not \in W$,
$e_t:W^{\perp} \rightarrow V^*$
has rank exactly $i-2$.   
Also $M(t) \not \in e_t(\barW^\perp)$, since that would imply
$M(t)(v) =0$. Thus
$e_t(Z) \subset V^*$ has dimension $i-1$
and so is equal to $e_t(X)$.
Now let $x$ be an arbitrary element of $X$. Since
$e_t(X)=e_t(Z)$, there exists an element $z \in Z \subset W^\perp$ such
that $x(t)=z(t)$.  But then, for any $w \in W$, 
$$
x(w)(t)=x(t)(w)=z(t)(w)=z(w)(t)=0. \leqno(23)
$$
Since $t$ was chosen
arbitrarily from an open set, $x(w)$ must be zero.
Since $x$ and $w$ were arbitrary, $X \subset W^\perp$. But $X$ and
$W^\perp$ have the same dimension, so $X=W^\perp$.
\qed
\enddemo

\subhead \S 5 The tangent space is constant \endsubhead
\proclaim{5.1 Theorem} Let $X \subset \hg$ be
a variety of dimension $\frac{i(i-1)}{2}$, with $i \geq 3$, 
such that at every smooth point $A \in X$ there is 
a $(g-i+1)$-plane $W(A) \subset \bC^g$ such that
$$
T_A X =W(A)^{\perp} \subset T_A \hg = S_g. \leqno(24)
$$ 
Then $W(A) = W$ is constant and $X \subset \hg$
is the affine subspace 
$$
X = \{M \in \hg| M|_W = \tau|_W\} \leqno(25)
$$
for some fixed $\tau \in \hg$.
\endproclaim

\demo{Proof} Once $W(A)$ is constant, the tangent
space is constant and equation (25) follows.
Let $I:= \{1,\dots,i-1\}$, $K = \{i,\dots,g\}$. 
Pick a smooth point $A_0$. After changing basis we can assume
$W(A_0)$ is the span
of $\{e_k\}, k\in K$, where $e_1, \ldots, e_g \in \bC^g$ are
the standard basis elements. Then 
on an open set $U$ around $A$ 
the coordinates $A_{ab}$ of the
upper left $(i-1) \times (i-1)$ block are analytic
coordinates for $X$, and there exist holomorphic
functions $z_{a,t}, a \in I, t \in K$ so that 
$$
v_t(A) := e_t + \sum_{s \in I} z_{s,t} e_s  \leqno(26)
$$
form a basis of $W(A)$. 
We will show that $\partial_{ab} z_{s,t} =0$ for all
$a,b,s \in I, t \in K$ 
(where $\partial_{ab}$ means
differentiation with respect to the coordinate $A_{a,b}$).
At $A \in U$ 
a basis for $T_A X \subset S_g$ is given by the
matricies $\partial_{ab} A$. 
Note that the upper-left-hand entries of 
$\partial_{ab} A$ are particularly simple. If $a,b,c,d \in I$, then
$$ (\partial_{ab} A)_{cd} = \hbox{1 if ($a=c$ and $b=d$) 
or ($a=d$ and $b=c$), 0 otherwise}. \leqno(27)
$$
By assumption
$(\partial_{ab} A) \cdot v_t = 0$. Looking at
the first $i$ entries of this product gives:
$$
\align
&\partial_{ab} A_{c,t} =0 \quad a,b,c \in I, t \in K, c \not \in
\{a,b\} \\
&z_{a,t} + \partial_{ab} A_{b,t} = 0 \quad a,b \in I, t \in K \tag{28}
\endalign 
$$

If $i > 3$, or if $i =3$ and $a=b$, then 
there exists 
$c \in I \setminus \{a,b\}$. 
Thus for $s\in I, t \in K$,
$$ \partial_{ab} z_{s,t} = -\partial_{ab} (\partial_{sc}
A_{c,t})
= - \partial_{sc} (\partial_{ab} A_{c,t})=0. \leqno(29)
$$

Finally, if $i=3$ and $a,b \in I$ we have 
$$
\align
\partial_{ab} z_{a,t} &= -\partial_{ab}(\partial_{aa} 
A_{a,t}) \\
                 &= -\partial_{aa}(\partial_{ba} A_{a,t}) \\
                 &= \partial_{aa} z_{b,t} = 0  \qed  \tag{30} 
\endalign 
$$
\enddemo

\subhead \S 6 Compactness \endsubhead
Here we argue that the image of $X(W,\tau)$ 
in $\ag$ cannot be compact. The meaning
of $X(W,\tau)$ is clearer in an alternative
realization of $\hg$: 

$\hg$, with its $\Sp(2g,\bR)$ action
is canonically identified with the 
set of complex structures on 
$\bR^{2g} = \bC^g$ compatible with the
standard symplectic form $(, )$, where
$\Sp(2g,\bR)$ acts on complex structures  by conjugation.
For $M \in \hg$ the corresponding complex
structure on $\bR^{2g}$ is given by identifying
$\bR^{2g}$ with $\bC^g$ by sending
a column vector $x \in \bC^g$ to
$$
f_{M}(x) = (\re(x),-\re(M x)). \leqno(31)
$$
The above, and the next Lemma, are all immediate
from \cite{Mumford83,4.1}

\proclaim{6.1 Definition-Lemma} For
$M \in X(W,\tau)$
$$
V:= f_M(W) = f_{\tau}(W) \leqno(32)
$$
is a non-degenerate subspace of $\bR^{2g}$
with a fixed complex structure compatible
with the restriction of $(,)$. 
Under the
above realization of $\hg$, 
$X(W,\tau)$ is identified with the set of complex
structures on $\bR^{2g}$, compatible
with $(,)$, extending the complex structure
structure on $V$.
This is naturally identified
with compatible complex structures on 
$(V^{\perp},(,)|_{V^{\perp}})$, and thus with
$\vhg r.$, where $r$ is the complex codimension
of $W$. 
\endproclaim

Let $H$ (resp $\Gamma \subset H$) be the subgroup of
$\Sp(2g,\bR)$ 
(resp $\Sp(2g,\bZ)$)
that preserves $X(W,T)$ -- or equivalently by
the above, the subgroups that preserve, and act
complex linearly on $V$. Restriction to $V,V^{\perp}$
identifies $H$ with $\U(s) \times \Sp(2g,\bR)$,
and under this identification, $\U(s) \times \{1\}$
is the stabilizer of $\tau$. Recall
that a subgroup of a semi-simple Lie group
is called a {\bf lattice} if it is discrete,
and the quotient has finite volume. It is called
cocompact iff the quotient is compact.
The following is clear:

\proclaim{6.2 Lemma} Assume the image, $Z$,
of $X(W,\tau)$ in
$\ag$ is a closed analytic subvariety. Then 
$$
\Gamma \backslash (H / \U(s)) =  \Gamma \backslash  X(W,\tau) \leqno(33)
$$
is the normalization of $Z$,
$\Gamma \subset H$ is a lattice, and cocompact iff
$Z$ is compact.
\endproclaim

As $\U(s)$ is compact, $\Gamma \subset H$ is a 
lattice iff its image in $H/\U(s) = \Sp(2r,\bR)$
is a lattice. Furthermore, by the Borel Density
Theorem, \cite{Zimmer84,3.2.5}, a lattice
in a semi-simple Lie group without compact factors
is Zariski dense. 

Thus it sufficies to prove the following, which was done
jointly with (and mostly by) Scot Adams:

\proclaim{6.3 Theorem} Let $W \subset \bR^{2g}$ be
a real subspace, of codimension $2r$, together with
a complex structure compatible with $(,)|_W$.
Let $H \subset \Sp(2g,\bR)$ be the subgroup
of elements that preserves, and acts complex
linearly on, $W$. Let
$\Gamma := H \cap \Sp(2g,\bZ)$. The
image of the restriction map
$$
\Gamma \rightarrow \Sp(W^{\perp},\bR) = \Sp(2r,\bR) \leqno(34)
$$
is Zariski dense iff $W$ is spanned by integer vectors.
In this case the image is a non-cocompact lattice.
\endproclaim
\demo{Proof} 
Suppose first $W$ has a basis of integer vectors. 
Then the image of $\Gamma$ in 
$H /\U(s) = \Sp(2r,\bR)$ is commensurable with
(i.e. up to finite index subgroups equal to)
$\Sp(2r,\bZ)$, a non-cocompact lattice. 

Now suppose that $W$ is not spanned by integer vectors, and that
the image of $\Gamma$ is dense.
Let $G \subset H \subset \Sp(2g,\bR)$ be
the Zariski closure of $\Gamma$. 
Restriction  to $W^{\perp}$ 
gives 
a surjection 
$$
G(\bR) \rightarrow \Sp(W^{\perp},\bR) = \Sp(2r,\bR) \tag{35}
$$
which together with restriction to $W$  
defines
an embedding
$$
G(\bR) \subset \U(s) \times \Sp(2r,\bR). \tag{36}
$$

For any complex subspace $V \subset \bC^{2g}$, non-degenerate
with respect 
to the complex linear extension of $(, )$, let
$G_V \subset \Sp(2g,\bC)$ be the subgroup of
matricies that act trivially on $V$. Of course
restriction to $V^{\perp}$ identifies this
with $\Sp(V^{\perp},\bC)$. We indicate its
(complex) Lie algebra by 
$$
\sp_{V^{\perp}} \subset \sp(2g,\bC). \leqno(37)
$$
We abuse notation and refer to $W \otimes_{\bR} \bC$
as $W$ as well. 

Note that if a Lie subalgebra ${\goth h} \subset \sp(2g,\bC)$
is stable under $\Gal(\bC/L)$ for a subfield $L \subset \bC$,
then ${\goth h}$ is defined over $L$, i.e. is the extension
of scalars of a Lie subalgebra of $\sp(2g,L)$, which
we denote ${\goth h}(L)$. 
The subalgebra ${\goth h}(L) \subset {\goth h}$ is simply the subset
fixed by $\Gal(\bC/L)$. 

Since $W$ is not spanned by rational vectors, there exists
an element $\sigma \in \Gal(\bC/\bQ)$ for which $\sigma(W) \ne W$. Let
$V = \sigma(W)$.
Let $\lag \subset \sp(2g,\bC)$ be the complexified
Lie algebra of $G$.

By Levi's theorem, \cite{Serre65,4.1}, 
(35) has a section at the
level of real Lie algebras. Since there is no
non-trivial homomorphism $\sp(2r,\bR) \rightarrow \u(s)$
we conclude from (36) that 
$\sp_{W^{\perp}} \subset \lag$. Note the image
is an ideal, as it corresponds to the kernel
of the restriction
$$
G \rightarrow U(s) \subset \Sp(W,\bR). \leqno(38)
$$
If $V = \bar V.$, let 
$$
r = \sp_{V^{\perp}} + \sp_{W^{\perp}} \subset \sp(2g,\bC). \leqno(39)
$$
Since $\Gal(\bC/\bQ)$ preserves $\Gamma$, it
preserves $\lag$, and thus
$$
\sp_{V^{\perp}} = \sigma(\sp_{W^\perp})
\subset \lag \leqno(40)
$$
is another ideal. Since $\sp_{V^{\perp}}$
and $\sp_{W^{\perp}}$ are distinct simple
ideals, their intersection is trivial, so $r$ is the direct sum
$$
r = \sp_{V^{\perp}} \oplus \sp_{W^{\perp}} \subset \lag, \leqno(41)
$$
and is itself an ideal. 

If $V \neq \bar V.$, let
$$
r = \sp_{V^{\perp}} + \sp_{{\bar V.}^{\perp}} 
 + \sp_{W^{\perp}}
= \sp_{V^{\perp}} + {\bar {\sp_{V^{\perp}}}.} +
\sp_{W^{\perp}}  \subset \sp(2g,\bC). \leqno(42)
$$
By the same reasoning as before, each factor is a simple ideal,
so $r$ is the direct sum of its factors, and is an ideal. 

In either case, we have an ideal
$r \subset \lag$.
Since $r,\lag$ are preserved by complex conjugation,
we have an induced real ideal $r(\bR) \subset \lag(\bR)$.
By construction $r(\bR)$ is isomorphic to
$$
\align
\sp(2r,\bR) & \oplus \sp(2r,\bR) \text{ if } V = \bar V. \\
\sp(2r,\bC)  & \oplus \sp(2r,\bR) \text{ if } V \neq \bar V. \tag{43}
\endalign 
$$
By (35) there
is an induced surjection of real lie algebras
$$
r(\bR) \rightarrow \sp(2r,\bR) \leqno(44)
$$
whose kernel, by (36), is a subalgebra of
$\u(s)$. As $r(\bR)$ is semi-simple, the kernel is a direct sum of simple
factors of $r(\bR)$. Since no factor has a non-trivial
homomorphism to $\u(s)$,  we conclude $r$ has only
a single factor, which is a contradiction. \qed \enddemo

\subhead \S 7 Remarks on characteristic $p > 0$ \endsubhead

In positive characteristic the only
known codim $g$ complete subvariety of $\ag$
is the locus $Z \subset \ag \otimes {\Bbb F}_p$
of Abelian varieties of $p$-rank zero, discovered
by Oort. Results
of \cite{Koblitz75} imply that the tangent
space to $Z$ is 
just as in Corollary 2.6 (with $i=g$): 

\proclaim{7.1 Proposition (Koblitz)} 
At points $[A]$ in a Zariski dense subset
of $Z$ the following hold: The $p$-linear
Frobenius 
$$
H: H^1(A,\ring A.) \rightarrow H^1(A,\ring A.) \leqno(45)
$$
has one dimensional cokernel. Let 
$L_{[A]} \subset H^0(A,\Omega^1_A) = \bE_{[A]}$
be the dual line. $Z$ is (orbifold) smooth
at $[A]$ and its tangent space is given by
$$
T_{[A]} Z = L_{[A]}^{\perp} \subset S(\bE_{[A]},\bE^*_{[A]})
= T_{[A]} \ag. \leqno(46)
$$
\endproclaim
\demo{Proof} This follows from \cite{Koblitz75}. Here
we give a few details for the readers convenience. The
$p$-linear Frobenius gives a commutative diagram
$$
\CD
0 @>>> H^0(A,\Omega^1_A) @>>> H^1_{\dr}(A) @>>> H^1(A,\ring A.)
@>>> 0 \\
@. @V{0}VV @V{F}VV @V{H}VV @. \\
0 @>>> H^0(A,\Omega^1_A) @>>> H^1_{\dr}(A) @>>> H^1(A.\ring A.)
@>>> 0
\endCD \tag{47}
$$
The snake lemma induces a canonical map
$B: \ker(H) \rightarrow H^0(A,\Omega^1_A)$. This is
the (restriction to $\ker(H)$ of the)
map given by the matrix $B$ in \cite{Koblitz75}. The
argument on pages 188-189 shows that the tangent
space to $Z$ is the perpendicular to the (one dimensional) image 
of $B$. So it remains to show that the image of $B$ is
dual to the cokernel of $H$. This follows from the 
fact that the image of $F$ is isotropic for the 
canonical pairing on $H^1_{\dr}$. \qed \enddemo

Given the parallel between Theorem 7.1 and our argument
in characteristic zero, it is natural to wonder:

\proclaim{7.2 Question} If $Z \subset \ag \otimes {\Bbb F}_p$
is complete and of codimension $g$ is the tangent space to
$Z$ as in Corollary 2.6? Is Oort's example the only possibility?
\endproclaim

\subhead \S 8 $\mgc$ \endsubhead
\demo{Proof of 1.2.1 } By \cite{Diaz87} a 
compact subvariety of $\mgc$ or $\mgco$ has
codimension at least $g$, for any $g$, and
a compact subvariety of $\omg$ or $\omgo$ has
codimension at least $2g-1$. Now assume
$g \geq 3$ and let $Z \subset \mgc$ be
a compact subvariety of codimension at
most $g$. 
We have a regular map
$$
\mgc \rightarrow \ag
$$
with zero dimensional fibres outside of 
$\partial \mgc$ (meaning the complement 
of $\omg \subset \mgc$). 
Suppose first $g=3$. By (1.2)
$Z \subset \partial M_3^c$ and so 
projects to a complete
surface in $M_2^c$, violating Diaz's bound.
So we may assume $g \geq 4$. We proceed
by induction.

By Diaz's bounds $Z$ must meet the 
boundary. Let $Z_i$ be the intersection
of $Z$ with $\delta_i$.

If $Z_i$ is non-empty we have 
by Diaz's bounds 
$$
\align
g  &\geq \codim(Z,\mgc) \geq \codim(Z_i,\delta_i) \\
&\geq \codim(\pi_i(Z_i), M_{i,1}^c) + 
\codim(\pi_{g-i}(Z_i), M_{g-i,1}^c) 
\endalign
$$
where $\pi_k$ indicates the projection onto
the factor $M_{k,1}^c$. Furthermore, by
induction, we can replace the last term by
$g+1$ if either $i$ or $g-i$ is at least $3$. 
Thus the only possibility is that 
$g=4$ and $Z$ meets only $\delta_2$. But then
its image in $M_{2,1}^c$ is a complete surface
contained in $M_{2,1}$, contradicting Diaz's
bounds. \qed \enddemo


\refstyle{C}

\end
\Refs
\widestnumber\key{Griffiths84}
\ref \key {BD85} 
      \by T. ~Brocker and T. tom Dieck
      \book Representations of Compact Lie Groups
      \publ Springer
      \yr 1985
\endref

\ref \key {BT82} 
     \by  R. ~Bott and L. ~Tu
     \book Differential Forms in Algebraic Topology
     \publ Springer
     \yr 1982
\endref
\ref \key{CP90} \by E. ~Colombo and G. P ~ Pirola
     \paper Some density results for curves with nonsimple Jacobians
      \jour Math. Ann.
      \vol 288
      \yr 1990
      \pages 161--178
\endref
\ref \key {Diaz87} \by S. ~Diaz
     \paper Complete subvarieties of the moduli space of smooth curves
     \yr 1987
     \inbook Algebraic geometry, Bowdoin, 1985 (Brunswick, Maine, 1985)
     \pages 77--81
     \publ Amer. Math. Soc.
      \publaddr Providence, RI
\endref

\ref \key {EV02} \by H. ~Esnault and E. Viehweg
     \paper Chern classes of Gauss-Manin bundles of weight
$1$ vanish
      \paperinfo preprint math.AG/020103
      \yr 2002
\endref
\ref \key {Faber99} \by  C. ~Faber
     \paper A Conjectural Description of the Tautological
Ring of the Moduli Space of Curves
     \inbook Moduli of Curves and Abelian Varieties
      \pages 109--125
      \yr 1999
\endref   
\ref \key {FL99} \by  C. ~Faber and E. ~Looijenga
     \paper Remarks on Moduli of Curves
     \inbook Moduli of Curves and Abelian Varieties
      \pages 23--39
      \yr 1999
\endref    

\ref\key {FP00} \by C.~Faber and R. Pandharipande
     \paper Logarithmic Series and Hodge Integrals in the
Tautological Ring
     \jour Michigan Math. Jour.
     \vol 48
     \yr 2000
     \pages 21--239
\endref


\ref \key {Fulton84} \by W. ~Fulton
      \book Intersection Theory
      \publ Springer-Verlag
      \yr 1984
\endref


\ref \key {Griffiths84} \by P. ~Griffiths
      \paper Curvature Properties of the Hodge bundles
      \inbook Topics in transcendental algebraic geometry 
      \jour Ann. of Math. Stud.
      \vol 106
      \pages 29--49, 
      \yr 1984
\endref
\ref \key {GH78} \by  P. ~Griffiths and J. ~Harris
       \book Principles of Algebraic Geometry
        \publ Wiley-Interscience
        \yr 1978
\endref
\ref \key {G99} \by  G. ~van der Geer
      \paper Cycles on the Moduli Space of Abelian Varieties
      \inbook Moduli of Curves and Abelian Varieties
      \pages 65--88
      \yr 1999
\endref
\ref \key {GO99} \by  G. ~van der Geer
and F. ~Oort
     \paper Moduli of Abelian Varieties: A Short Introduction
and Survey
     \inbook Moduli of Curves and Abelian Varieties
      \pages 1--17
      \yr 1999
\endref    
\ref \key {Izadi98} \by E. ~Izadi
     \paper Density and completeness of subvarieties of 
moduli spaces of curves or abelian varieties
      \jour Math. Ann.
      \yr 1998
      \vol 310
      \pages 221-233
\endref

\ref \key {Mumford83} \by  D. ~Mumford
      \paper Towards an enumerative geometry of the moduli
space of curves
      \inbook Arithmetic and Geometry vol. II
      \pages 271-328
      \jour Progress in Math.
      \vol 36
      \yr 1983
\endref
\ref\key {Kempf91} \by G. ~Kempf
    \book Complex Abelian Varieties and Theta Functions
    \publ Springer-Verlag
    \yr 1991
\endref
\ref \key {Koblitz75} \by N. ~Koblitz
     \paper $p$-adic variation of the zeta-function over 
families of varieties defined over finite fields
     \jour Compos. Math.
     \vol 31
     \yr 1975
     \pages 119--218
\endref
\ref \key {Serre65} \by J.P. ~Serre
      \book Lie algebras and Lie groups
      \bookinfo Lectures given at Harvard University, 1964
      \publ W.A. Benjamin, Inc.
      \yr 1965
\endref
\ref \key {Zimmer84} \by R. ~ Zimmer
     \book Ergodic theory and semisimple groups
      \publ  Birkhauser Verlag
      \yr 1984
\endref
\endRefs
